\documentclass[onecolumn,11pt]{article}

\usepackage{algorithm}
\usepackage{algorithmic}
\usepackage{multirow,bigdelim}
\usepackage{amsmath}
\usepackage{bm}
\usepackage{amssymb}
\usepackage{amscd}
\usepackage[dvipdfm]{graphicx,color}
\usepackage{amsfonts} 
\usepackage{fullpage} 
\usepackage[tableposition=top]{caption}
\newtheorem{define}{Definition}[section]
\newtheorem{thm}{Theorem}
\newtheorem{lem}{Lemma}

\newtheorem{prop}[define]{Proposition}

\makeatletter
   
   \@addtoreset{table}{section}
  \makeatother
\makeatletter
    
    \@addtoreset{figure}{section}
  \makeatother

\title{Perturbation bounds of eigenvalues of Hermitian matrices with block structures}

\author{ Yuji Nakatsukasa \\\\Department of Mathematics, University of California, Davis\\ynakatsukasa@ucdavis.edu} 
\date{}
\begin{document}
\large
\maketitle
\normalsize
\begin{abstract} 
We derive new perturbation bounds for eigenvalues of Hermitian matrices with block structures. The structures we consider range from a standard 2-by-2 block form to block tridiagonal and tridigaonal forms. The main idea is the observation that an eigenvalue is insensitive to componentwise perturbations if the corresponding eigenvector components are small. We show that the same idea can be used to explain two well-known phenomena, one concerning extremal eigenvalues of Wilkinson's matrices and another concerning the efficiency of aggressive early deflation applied to the symmetric tridiagonal QR algorithm. 
\end{abstract}
\section{Introduction}
Theory of eigenvalues of Hermitian matrices is a well-studied subject, with many aesthetically pleasing results available, such as the max-min characterization, Cauchy's interlacing theorem and Weyl's theorem \cite[Ch.4]{stewart-sun:1990} \cite[Ch.8]{Golubbook}, \cite[Ch.4]{demmelbook}. 
Here we are concerned with eigenvalue perturbation bounds, and first note that for an eigenvalue $\lambda$ of a given unstructured Hermitian matrix $A$ and a general $E$ such that $\|E\|_2$ is bounded by a known constant, Weyl's theorem gives the best possible bound that is attainable.

The goal of this paper is to specialize in block Hermitian matrices and derive perturbation bounds of the eigenvalues that are sharper than general bounds, such as Weyl's theorem.

Much work has been done in this direction as well. For example,  a well-known quadratic residual bound exists  \cite{parlettsym,mathiasquad} that relates the eigenvalues of two 2-by-2 block Hermitian matrices of the form
\begin{equation}  \label{origin}
A=  \begin{bmatrix}A_{1}&0\\0&A_{2}\end{bmatrix}  \quad \mbox{and}\quad \widehat A=  \begin{bmatrix}A_{1}&E^H\\E&A_{2}\end{bmatrix}, \end{equation}
if $\lambda_i(A)\notin \lambda(A_1)$ is the $i$th eigenvalue of $A$ then 
\begin{align}
|\lambda_i(A)-{\lambda_i}(\widehat{A})|&\leq \frac{\|E\|_2^2}{\min_j|\lambda_i(A)-\lambda_j(A_2)|}.\label{standquadvar}
\end{align}
Here, $\lambda_i(X)$ denotes the $i$th smallest eigenvalue of a Hermitian matrix $X$ and  $\lambda(X)$ denotes the set of $X$'s eigenvalues. 

When $\|E\|_2$ is small, \eqref{standquadvar} provides much tighter 
bounds than linear eigenvalue perturbation bounds do, such as the well-known Weyl's theorem \cite{demmelbook}, which gives  $|\lambda_i(A)-\lambda_i(\widehat A)|\leq \|E\|_2$. We note that \cite{rcli05} gives an improved bound that is always sharper than both the bound by Weyl's theorem and \eqref{standquadvar}. 

The contribution of this paper is that  we present a new framework for deriving new eigenvalue perturbation bounds, which is to first obtain bounds for the relevant  eigenvector components and use them to arrive at bounds of eigenvalues. 
Based on this framework we give new bounds for the 2-by-2 block form as in \eqref{origin}, but unlike \eqref{origin} we assume no zero submatrix. Specifically, we study the difference between eigenvalues of $A$ and $A+E$, where 
\begin{equation}  \label{origin2}
A=  \begin{bmatrix}A_{11}&A_{21}^T\\A_{21}&A_{22}\end{bmatrix}  \quad \mbox{and}\quad E=  \begin{bmatrix}E_{11}&E_{21}^H\\E_{21}&E_{22}\end{bmatrix}, \end{equation}
We then turn to the more specialized and  well-studied tridiagonal cawse and show a tight bound results from our idea when a target eigenvalue is disjoint from many Gerschgorin disks.
To demonstrate the sharpness of our approach, we show that our framework successfully explains the following two well-known phenomena: (i) Wilkinson's matrices have many pairs of nearly equal eigenvalues. (ii) Aggressive early deflation applied to the symmetric tridiagonal QR algorithm deflates many eigenvalues even when no subdiagonal element is negligibly small. 


The structure of this paper is as follows. In section \ref{prelim} we survey and derive some necessary results, and present our basic idea . Section \ref{secbounds} treats the 2-by-2 block case and proves a new bound that is sometimes tighter than any known bound. 
Section \ref{tridiagcase} deals with the tridiagonal case, in which we investigate the above two case studies. 

 Notations: $\lambda_i(X)$  denotes the $i$th smallest eigenvalue of a Hermitian matrix $X$. For simplicity  we use $\lambda_i,\lambda_i(t)$ and $\widehat\lambda_i$ to denote the $i$th eigenvalue of $A,A+tE$ and $A+E$ for $t\in[0,1]$ and $A$ and $E$ as defined in \eqref{origin2} respectively.
  $\lambda(A)$ denotes the set of the eigenvalues of a Hermitian matrix $A$. $\sigma_i(B)$ denotes the $i$th largest singular value of a general matrix $B$. 
We only  use the spectral norm $\|\cdot \|_2$. 

\section{Preliminaries}\label{prelim}
In this section we review some results that will be necessary for our analysis. 
 We first recall the partial  derivative of simple eigenvalues \cite{stewart-sun:1990}.
\begin{lem}\label{lem0}
Suppose $A$ and $E$ are Hermitian as in \eqref{origin2}. 
Denote by $\lambda_i(t)$ the $i$th eigenvalue of $A+tE$ such that $(A+tE)x^{(i)}(t)=\lambda_i(t)x^{(i)}(t)$ where $\|x^{(i)}(t)\|_2=1$ for some $t\in [0,1]$. If $\lambda_i(t)$ is simple, then\footnote{For simplicity hereafter we discard the superscript $(i)$ in the eigenvector $x^{(i)}(t)$.}
\begin{equation}\label{pert}
\frac{\partial \lambda_i(t)}{\partial t}= x(t)^HE x(t).\end{equation}
\end{lem}
Since $\lambda_i(0)=\lambda_i(A)=\lambda_i$ and $\lambda_i(1)=\lambda_i(A+E)=\widehat\lambda_i$, from \eqref{pert} it follows that if $\lambda_i(t)$ is simple for all $0\leq t\leq 1$, then 
\begin{align}
|\lambda_i -\widehat\lambda_i|&=\left|\int_0^1 x(t)^HE x(t)dt\right|  \label{int}\\
&\leq \left|\int_0^1 x_1(t)^HE_{11} x_1(t)dt\right|+2\left|\int_0^1 x_2(t)^HE_{21} x_1(t)dt\right|+\left|\int_0^1 x_2(t)^HE_{22} x_1(t)dt\right|,  \label{int2}
\end{align}
where we partitioned $x(t)=\begin{bmatrix}x_1(t)\\x_2(t)\end{bmatrix}$ for $x_2(t)\in\mathbb{C}^{k}$.
The key observation here 
 is that the latter two terms in \eqref{int2} are small 
 if $\|x_2(t)\|_2$ is small for all $0\leq t\leq 1$. 
The next Lemma gives a useful upper bound for $\|x_2(t)\|_2$. 


\begin{lem}\label{lem}
Suppose that $\lambda_i\notin \lambda(A_{22})$ is the $i$th eigenvalue of $A$ as defined in \eqref{origin2}\footnote{If $\lambda\notin \lambda(A_{11})$ we simply swap the subscripts $1$ and $2$ in the following arguments.}. 
Let $Ax=\lambda_i x$  such that $\|x\|_2=1$. 
Then, denoting  $x=\begin{bmatrix}x_1\\x_2\end{bmatrix}$ for $x_2\in\mathbb{C}^{k}$ we have
\begin{equation}  \label{lem1}
\|x_2\|_2\leq \frac{\|A_{21}\|_2}{\min_{j}|\lambda_i-\lambda_j(A_{22})|}. 
\end{equation}
\end{lem}
{\sc proof.} The bottom $k$ rows of $Ax=\lambda x$ is 
\[A_{21}x_1+A_{22}x_2=\lambda_i x_2,\]
so we have 
\[x_2=(\lambda_i I-A_{22} )^{-1}A_{21}x_1.\]
Taking norms we get 
\[\|x_2\|_2\leq \|(\lambda_i I-A_{22} )^{-1}\|_2\|A_{21}\|_2\|x_1\|_2\leq \frac{\|A_{21}\|_2}{\min_{j}|\lambda_i-\lambda_j(A_{22})|}.\ (\because \|x_1\|_2\leq \|x\|_2=1)\]
 \hfill $\square$
\bigskip

We note that \eqref{lem1} is valid for any $\lambda_i$ and its eigenvector $x$, whether or not $\lambda_i$ is a  multiple eigenvalue. Thus for the multiple case, all the vectors that span the corresponding eigenspace satisfies \eqref{lem1}.  

However, Lemma \ref{lem0} assumes that $\lambda_i$ is a simple eigenvalue of $A$ to derive the partial derivative of $\lambda_i$ with respect to $t$. 
 Special treatment is needed to get the derivative of multiple eigenvalues, and this is shown in the appendix. It turns out that 
everything that we discuss carries over, 
in that when $\lambda_i(t)$ is multiple, \eqref{pert} still holds for a certain eigenvector $x(t)$ of $\lambda_i(t)$. We defer the treatment of multiple eigenvalues to the appendix, because they only cause complications to the analysis that are not fundamental to the eigenvalue behavior. 
Hence for simplicity  we assume that $\lambda_i(t)$ is simple for all $t$, so that the normalized eigenvector is unique (up to a factor $e^{i\theta}$). 
\section{2-by-2 block case}\label{secbounds}
In this section we derive refined eigenvalue perturbation bounds by combining Lemmas \ref{lem0} and \ref{lem}. 

Consider the matrix $A+tE$ and its $i$th eigenvalue $\lambda_i(t)$ such that $(A+tE)x(t)=\lambda_i(t)x(t)$. 
The key observation is that if $\min_{j}|\lambda_i-\lambda_j(A_{22})|>\|E_{22}\|_2$ then \eqref{lem1} provides an upper bound for $\|x_2 (t)\|_2$ 
 for all $t\in[0,1]$: 
\begin{equation}
  \label{bound1}
\|x_2 (t)\|_2\leq \frac{\|A_{21}\|_2+\|E_{21}\|_2}{\min_{j}|\lambda_i-\lambda_j(A_{22})|-2\|E\|_2}\ (\equiv \tau_i),\quad 0\leq t\leq 1.
\end{equation}
This is verified simply by letting $A\leftarrow A+tE$ in \eqref{lem1}, which gives 
\[\|x_2(t)\|_2\leq \frac{\|A_{21}+tE_{21}\|_2}{\min_{j}|\lambda_i(t)-\lambda_j(A_{22}+tE_{22})|}, 
\]
and using the facts that 
$\|A_{21}+tE_{21}\|_2\leq \|A_{21}\|_2+t\|E_{21}\|_2$  and $\min_{j}|\lambda_i(t)-\lambda_j(A_{22}+tE_{22})|\geq  \min_{j}|\lambda_i(0)-\lambda_j(A_{22})|-\|E\|_2-\|E_{22}\|_2$ for all $0\leq t\leq 1$, because $|\lambda_i(t)-\lambda_i(0)|\leq t\|E\|_2$  by Weyl's theorem. 

\bigskip

Now we present our main result. 
\begin{thm}\label{main}
Let  $\lambda_i,\widehat\lambda_i$ be the $i$th eigenvalue of $A$ and $A+E$ as in \eqref{origin2} respectively, and define $\displaystyle\tau_i= \frac{\|A_{21}\|_2+\|E_{21}\|_2}{\min_{j}|\lambda_i-\lambda_j(A_{22})|-2\|E\|_2}    $ as in \eqref{bound1}. For each $i$, if $\tau_i>0$ then
 \begin{align}
\left|\lambda_i -\widehat\lambda_i\right|&\leq \|E_{11}\|_2+2\|E_{21}\|\tau_i +\|E_{22}\|_2\tau^2_i  \label{result}.
  \end{align}
\end{thm}
{\sc proof.}
Substituting \eqref{bound1} into \eqref{int2} we get
\begin{align*}
\left|\lambda_i -\widehat\lambda_i\right|
&\leq\left|\int_0^1 \|E_{11}\|_2\|x_1 (t)\|_2^2dt \right|+2\left|\int_0^1\|E_{21}\|_2 \|x_1 (t)\|_2 \|x_2 (t)\|_2dt \right|+\left|\int_0^1 \|E_{22}\|_2 \|x_2 (t)\|_2^2dt\right|\\
&\leq \|E_{11}\|_2+2\|E_{21}\|_2\tau_i +\|E_{22}\|_2\tau_i^2,
\end{align*}
which is \eqref{result}.\hfill $\square$
\bigskip

Two remarks are in order regarding the theorem. 
\begin{itemize}
\item Theorem \ref{main} is tighter than the Weyl bound $\|E\|_2$ only when $\tau_i<1$, which is $\min_{j}|\lambda_i-\lambda_j(A_{22})|>\|E_{22}\|_2+\|A_{21}\|_2+\|E_{21}\|_2$. If $\lambda$ is far from the spectrum of $A_{22}$ and $\|E_{11}\|_2$ is small, then Theorem \ref{main} is much tighter than Weyl's theorem. 

\item When $A_{21}$ and $E_{11}$ are zero but $E_{22}$ and $E_{21}$ are nonzero,  \eqref{result} reveals  that $\lambda_i$ is particularly insensitive to the perturbation $E_{22}$: in \eqref{result}, the term involving  $\|E_{22}\|_2$ becomes proportional to $\|E_{22}\|_2\|E_{21}\|_2^2$, which scales \emph{cubically} with $\|E\|_2$. 

For example, consider the $n$-by-$n$ matrices $\begin{bmatrix}A_{11}& \delta\\ \delta& \epsilon  \end{bmatrix}$ and $\begin{bmatrix}A_{11}& \delta\\ \delta& 0  \end{bmatrix}$ where $A_{11}$ is nonsingular. 
These matrices have $n-1$ eigenvalues that match up to $O(\epsilon\delta^2)$, and one eigenvalue that matches up to $\epsilon$. 
\end{itemize}

\section{Tridiagonal case}\label{tridiagcase}
We now turn to the symmetric tridiagonal case, and suppose that 
\begin{equation}  \label{tridef}
A=\begin{bmatrix}a_1&b_1&&\\b_1&\ddots&\ddots&  \\ &\ddots&\ddots&b_{n-1}\\&&b_{n-1}&a_{n}\end{bmatrix}\quad\mbox{and}\quad 
E=\begin{bmatrix}e_1&f_1&&\\f_1&\ddots&\ddots&  \\ &\ddots&\ddots&f_{n-1}\\&&f_{n-1}&e_{n}\end{bmatrix}, 
\end{equation}
where we assume without loss of generality that $b_i>0$ for all $i$. 
Our aim is to show that we can derive tighter results based on the same idea as in the previous section. We demonstrate the sharpness of our approach by considering the two case studies:
\begin{enumerate}
\item Explain why two largest eigenvalues of Wilkinson's matrices are nearly equal.
\item Explain why aggressive early deflation can deflate eigenvalues as ``converged'' when applied to the symmetric tridiagonal QR algorithm. 
\end{enumerate}
\subsection{Basic idea}
When  $A$ and $E$ are symmetric tridiagonal as defined in \eqref{tridef}, one can of course use Theorem \ref{main} to get eigenvalue perturbation bounds. However the tridiagonal structure enables us to refine the bound in Lemma \ref{lem}, which in turn yields tighter  eigenvalues bounds. 

Specifically, for symmetric tridiagonal $A$, the $k$th row of $Ax=\lambda x$ is 
\[b_{k-1}x_{k-1}+a_kx_k+b_{k}x_{k+1}=\lambda x_k,\]
so we have
\[|\lambda-a_k||x_k|= \left|b_{k-1}x_{k-1}+b_{k}x_{k+1}\right|\leq (|b_{k-1}|+|b_{k}|)\max(|x_{k-1}|,|x_{k+1}|).\]
Therefore if $|\lambda-a_k|\geq |b_{k-1}|+|b_{k}|$ then $|x_k|\leq \max(|x_{k-1}|,|x_{k+1}|)$ and 
\[|x_k|\leq \frac{|b_{k-1}|+|b_{k}|}{|\lambda-a_k|}\max(|x_{k-1}|,|x_{k+1}|).\]
Note that the condition $|\lambda-a_k|\geq |b_{k-1}|+|b_{k}|$ implies that  the disjointness of a Gerschgorin disk from an eigenvalue implies that the eigenvector components are decaying. This observation  was made in \cite{parlettger}. 

Below we show how this idea can be used in practice. 
\subsection{Eigenvalues of Wilkinson's matrix}\label{wilk}
The well-known Wilkinson's matrix \cite{wilkinson:1965}, whose famous $2n+1=21$ case is
\begin{equation}  \label{wilkmat}  
W_{21}^+=  \begin{bmatrix}10&1&&&&&\\1&9&\ddots&&&&\\&\ddots&\ddots&1&&&&&\\&&1&1&1&&&&\\&&&1&0&1&&&\\&&&&1&1&\ddots&&\\&&&&&\ddots&\ddots&\ddots&\\&&&&&&\ddots&9&1\\&&&&&&&1&10  \end{bmatrix}.
\end{equation}
 Such matrices are known to have many pairs of extremely close eigenvalues: for example, the two largest eigenvalues agree up to about $7\times 10^{-14}$. \cite[p.308]{wilkinson:1965} notes that in general the two largest eigenvalues of the matrix $W_{2n+1}^+$ agree up to roughly  $(n!)^{-2}$, but does not explain this in detail. 
We shall give an explanation using the ideas we described in this paper. Define $(2n+1)$-by-$(2n+1)$ matrices $A$ and $E$ such that $A+E=W_{21}^+$ by 
\begin{equation}  \label{wilkmatdiv}  
A=  \begin{bmatrix}10&1&&&&&\\1&9&\ddots&&&&\\&\ddots&\ddots&1&&&&&\\&&1&1&0&&&&\\&&&0&0&0&&&\\&&&&0&1&1&&\\&&&&&1&\ddots&\ddots&\\&&&&&&\ddots&9&1\\&&&&&&&1&10  \end{bmatrix},
E=
 \begin{bmatrix}&&&&&&\\&&\ddots&&&&\\&\ddots&\ddots&\ddots&&&&&\\&&\ddots &&1&&&&\\&&&1&0&1&&&\\&&&&1&&\ddots&&\\&&&&&\ddots&\ddots&\ddots&\\&&&&&&\ddots&&\\&&&&&&&&  \end{bmatrix}.
\end{equation}
Note that $A$ has $10$ (in general $n$) pairs of multiple eigenvalues with multiplicity 2, and a simple eigenvalue 0. 
We show that the large eigenvalues  of $A$ (those close to 10) are extremely insensitive to the perturbation  $E$, so that the two largest eigenvalues of $A+E$ must be very close to those of $A$, hence close to each other. 

Below we suppose $n>4$. First we consider the largest eigenvalue of  $A$, which we denote by $\lambda(>n)$. 
Define $x(t)=[x_1(t)\ x_2(t)\ \cdots\  x_{2n+1}(t)]$ such that $(A+tE)x(t)=\lambda(t) x(t)$ where  $\lambda(t)$ is a continuous function of $t$ with $\lambda(0)=\lambda$. 


Let us recall Lemma \ref{lem} and  consider refining the bound for $|x_{n}(t)|,|x_{n+1}(t)|$ and $|x_{n+2}(t)|$. 
Since $A$ is tridiagonal, from the $(n+1)$th row of  $(A+tE)x(t)=\lambda(t) x(t)$ we have 
\[\lambda(t)x_{n+1}(t)=t(x_{n}(t)+x_{n+2}(t)),\]
hence 
\begin{equation}  \label{xn}
|x_{n+1}(t)|\leq \frac{2t\max(|x_{n}(t)|,|x_{n+2}(t)|)}{|\lambda(t)|}
,\quad  t\in[0,1].  
\end{equation}
First we consider the case $|x_{n}(t)|>|x_{n+2}(t)|$, in which case  we also have $|x_{n}(t)|>|x_{n+1}(t)|$  in view of \eqref{xn}. 
From the $(n-1)$th row of $(A+tE)x(t)=\lambda(t) x(t)$ we similarly get 
\begin{equation}  \label{n-1}
|x_{n-1}(t)|\leq \frac{t(|x_{n-2}(t)|+|x_{n}(t)|)}{|\lambda(t)-1|},  \quad t\in[0,1].   
\end{equation}
Now since $n<\lambda(t)<n+1$ for all $t\in[0,1]$\footnote{We can get $n<\lambda(t)<n+1$ by first following the same argument using $n+2>\lambda(t)>n-1$.}  we must have $|x_{n-2}(t)|>|x_{n-1}(t)|>|x_{n}(t)|$. 
Substituting this into \eqref{n-1} yields $\displaystyle |x_{n-1}(t)|\leq \frac{t|x_{n-2}(t)|}{|\lambda(t)-1|-t}$. 
Therefore we have
%
\[|x_{n-1}(t)|\leq \frac{t|x_{n-2}(t)|}{n-2}\leq \frac{|x_{n-2}(t)|}{n-2},\quad t\in[0,1].\]
By a similar argument we find that 
 \begin{equation}   \label{gen}
|x_{n-i}(t)|\leq \frac{t|x_{n-i-1}(t)|}{n-i-1}\leq \frac{|x_{n-i-1}(t)|}{n-i-1} \quad \mbox{for}\quad 1\leq i\leq n-2,\  0\leq t\leq 1,
 \end{equation}
so together with \eqref{xn} we get
 \begin{equation}   \label{boundtight}
|x_{n+1}(t)|   \leq \frac{2t}{n}|x_2|\prod_{i=1}^{n-2} \frac{1}{n-i-1}  \leq \frac{2t}{n}\prod_{i=1}^{n-2} \frac{1}{n-i-1}
 \end{equation}
and
 \begin{equation}   \label{boundtight2}
|x_{n}(t)|\leq t\prod_{i=1}^{n-2} \frac{1}{n-i-1}. 
 \end{equation}
When $n=10$ we have $\delta_0<\delta_1<5\times 10^{-5}$. 
 which we now plug into \eqref{int2} 
 to get 
 \begin{align}
\left|\lambda_i(A+E) -\lambda_i(A)\right|&\leq  \left|\int_0^1 x (t)^HE x (t)dt\right|   \nonumber\\
&\leq  \int_0^1 2(|x_n (t)|+|x_{n+2} (t)|)|x_{n+1} (t)|dt   \nonumber\\
&\leq  \frac{4}{n}\left(\prod_{i=1}^{n-2} \frac{1}{n-i-1}\right)^2\int_0^1t^2dt\nonumber\\
&= \frac{4}{3n}\left(\prod_{i=1}^{n-2} \frac{1}{n-i-1}\right)^2.\label{fact}
 \end{align}
 The case $|x_{n-1}(t)|\leq |x_{n+1}(t)|$ can also be treated similarly, and we get the same result. 

We easily appreciate that the bound \eqref{fact} roughly scales as $1/n((n-2)!)^{2}$ as $n\rightarrow \infty$, which supports the claim in \cite{wilkinson:1965}.



We also note that by a similar argument we can prove that the $2\ell-1$th and $2\ell$th eigenvalues of $W_{2n+1}^+$ match to within
$1/(n-\ell+1)((n-\ell-1)!)^{2}$, which is small for small $\ell$, but not as small for larger $\ell$. 
Since this is an accurate description of what is well known about the eigenvalues of  $W_{2n+1}^+$, we conclude that this approach explains the observation that Wilkinson's matrix has many pairs of eigenvalues that are nearly equal.
\subsection{Aggressive early deflation applied to symmetric tridiagonal QR}\label{aggdef}
The aggressive early deflation strategy introduced in ~\cite{multi2} for the nonsymmetric Hessenberg QR algorithm, is known to greatly speed up the algorithm for computing the eigenvalues of a nonsymmetric matrix by deflating converged eigenvalues much earlier than a conventional deflation strategy does. Here we consider the symmetric tridiagonal QR, for which a similar (though perhaps not as dramatic) performance improvement is expected. 

The following is a brief description of aggressive early deflation applied to the symmetric tridiagonal QR. 
 Let  $A$ a tridiagonal matrix as defined in \eqref{tridef}. 
Denote by $A_{2}$ the  lower-right $k \times k$ submatrix of $A$, and let $A_{2}=VDV^T$  be an eigendecomposition, where the diagonals of $D$ are in decreasing order of magnitude. 
Then, we have
\begin{equation}
  \label{aggeq}
\begin{bmatrix}
I&\\&V
\end{bmatrix}^TA
\begin{bmatrix}I&\\&V
\end{bmatrix}=
\begin{bmatrix}\multicolumn{4}{c}{ \hbox{ \it{ \LARGE{ $A_1$} }} }&&t^T\\&&&t&&D\end{bmatrix}, 
\end{equation}
where $A_1$ is the  upper-left $(n-k) \times (n-k)$ submatrix of $A$, and
the vector $t$ is given by $t=b_{n-k}V(1,:)$ where
$V(1,:)$ is the first row of $V$. It often happens that many elements of $t$ are small, in which case aggressive early deflation regards  $D$'s corresponding eigenvalues as converged and deflate them. This is the case even when none of the subdiagonals of $A$ is particularly small. 

This means that many eigenvalues of the two matrices $A=\widehat A+E$ and $\widehat A=\begin{bmatrix}  A_{1}&0\\0&A_{2}\end{bmatrix}$ are almost equal, or equivalently that the perturbation of the eigenvalues by the $(n-k)$th subdiagonal $b_{n-k}$ is negligible. 
Here our aim is to explain why this is often the case. 
We do this by showing that under an assumption that is typically valid of a tridiagonal matrix appearing in the course of the QR algorithm, many eigenvalues of $A_2$ are perturbed only negligibly by $b_{n-k}$. Again we rely on our idea of bounding relevant eigenvector components. 

It is well-known that under mild assumptions the tridiagonal QR algorithm converges, in that the diagonals converge to the eigenvalues in descending order of magnitude, and the subdiagonal elements converge to
zero~\cite{trefbau}. In light of this, here we assume that the diagonals $a_i$ are roughly ordered in descending order of their magnitudes, and that the subdiagonals $b_i$ are small, so that for a target (small) eigenvalue $\lambda(A)$ of $A_2$, we have $|a_i-\lambda|>b_i+alpha$ for $1\leq i\leq n-k+j$ for some $j>0$. Here, $\alpha$ is a bound such that $|\lambda(t)-\lambda|\leq \alpha$ for all $0\leq t\leq 1$\footnote{We can safely let $\alpha=b_{n-k}$ which works, but we can get a much smaller bound for example by using the argument here recursively}. 

We bound the perturbation on $\lambda$ by $b_{n-s}$ by tracing $\lambda$  via \eqref{pert}. Specifically, defining
a continuous function $\lambda(t)$ of $t$ such that $(\widehat A+tE)\bm{x}(t)=\lambda(t)\bm{x}(t)$ where $\lambda(0)=\lambda$ and $\bm{x}(t)=[x_1(t)\ x_2(t)\ \cdots x_n(t)]$ is a unit vector for all $t$, 
we shall bound $|\lambda(1)-\lambda(0)|$. We shall prove the following. 
\begin{prop}\label{prop1}
Under the above notations and assumptions, 
\begin{equation}\label{pertsiki0}
|\lambda(1)-\lambda(0)|\leq \frac{b_{n-k}}{2}\eta_{n-k} \prod_{i=1}^{j}\eta_{n-k+i}^2,
\end{equation}  
\end{prop}
where $\displaystyle\eta_i=\frac{b_i}{|a_i-\lambda|-\alpha-b_i}$.\\
{\sc proof.} 
First note  from  \eqref{pert} that
\begin{equation}\label{pertsiki}
\lambda(1)-\lambda(0)=\int_0^{b_{n-k}}\bar x_{n-k}(t) x_{n-k+1}(t)dt.
\end{equation}
From the first row of $(\widehat A+tE)\bm{x}(t)=\lambda(t)\bm{x}(t)$ we have
\begin{align*}
\left(a_1-\lambda(t)\right)|x_1(t)|&=b_{1}|x_{2}(t)|,
\end{align*}
Hence 
\[ \frac{|x_{1}(t)|}{|x_2(t)|}=\frac{b_{1}}{a_1-\lambda(t)}<\frac{b_1}{|a_1-\lambda|-\alpha}<1,\]
 Next, assuming $|x_{i-1}(t)|<|x_{i}(t)|$ for some integer $i\leq  n-k+j$, from the $i$th row of  $(\widehat A+tE)\bm{x}(t)=\lambda(t)\bm{x}(t)$ we have
\begin{align*}
\left(a_k-\lambda(t)-b_{k-1}\frac{|x_{k-1}(t)|}{|x_k(t)|}\right)|x_k(t)|&=b_{k}|x_{k+1}(t)|,
\end{align*}
so it follows that 
\begin{align*} \frac{|x_{k}(t)|}{|x_{k+1}(t)|}&=\frac{b_{k}}{(a_{k}-\lambda(t)-b_{k-1}\frac{|x_{k-1}(t)|}{|x_k(t)|})}\\
&<\frac{b_k}{|a_{k}-\lambda|-\alpha-b_{k-1}}<1,
\end{align*}
where the last inequality follows from the assumption $|a_i-\lambda|>b_i+\alpha$. Hence, we have shown for $1\leq i\leq n-k+j$ that $|x_i(t)|\leq |x_{i+1}(t)|$ for all $0<t<1$.\\

 Therefore, we have
\[\frac{|x_{n-k}(t)|}{|x_{n-k+j}(t)|}\leq  \prod_{i=1}^{j}\frac{b_{n-k+i}}{|a_{n-k+i}-\lambda|-\alpha-b_{n-k+i-1}}.\]
 Since this holds for all $0<t<1$,  we conclude that 
\begin{align*}\lambda(1)-\lambda&=\lambda(1)-\lambda(0)=\int_0^{b_{n-k}}\bar x_{n-k}(t) x_{n-k+1}(t)dt\\
&\leq  \int_0^{b_{n-k}}\left|\frac{x_{n-k}(t)}{x_{n-k+j}(t)}\right| \left|\frac{x_{n-k+1}(t)}{{x_{n-k+j}(t)}}\right|dt\\
&\leq \frac{1}{2}\cdot\frac{b_{n-k}^2}{|a_{n-k}-\lambda|-\alpha-b_{n-k-1}}\cdot \left(\prod_{i=1}^{j}\frac{b_{n-k+i}}{|a_{n-k+i}-\lambda|-\alpha-b_{n-k+i-1}}\right)^2\\
&\leq \frac{b_{n-k}\eta_{n-k}}{2} \prod_{i=1}^{j}\eta_{n-k+i}^2.
\end{align*}
\hfill$\blacksquare$

Now let us analyze the result. The bound \eqref{pertsiki0} is a product of $j$ numbers $\eta_{n-k+i}^2$ for $1\leq i\leq j$, where each $\eta_{n-k+i}$ can be much smaller than $1$ if $A$ is a nearly diagonal matrix. Moreover, the above argument is valid as long as the  assumption on $\lambda$ that it is far from the $n-k+j$ diagonals of $A$ is true. 

\paragraph{Simple example} 
To illustrat the result, let $A$ be the $1000$-by-$1000$ tridiagonal matrix 
\begin{equation}\label{defa}
 A = \mbox{tridiag}\left\{
 \begin{array}{cccccccccccc}
      & 1 &     & 1 &    & . &   & 1 &   & 1 \\
100 &     & 999 &     & .  &   & . &     & 2 &     &1\\
      & 1 &     & 1 &    & . &   & 1 &   & 1 \\
 \end{array}
 \right\}.
\end{equation}
Let $k=100$, and let us target on an eigenvalue $\lambda$ of $A_{2}$ that is smaller than $10$ (there are at least 9 of them). For such $\lambda$ we can let $j=88$. 
Since $\displaystyle\eta_i=\frac{b_{n-k+i}}{|a_{n-k+i}-\lambda|-\alpha-b_{n-k+i-1}}\leq \frac{1}{|100-i-\lambda|-1}\leq \frac{1}{|90-i|-1}$ (where we used the safe bound $\alpha=1$), 
Proposition \ref{prop1} gives a bound
\begin{align*}
|\lambda(1)-\lambda(0)|&\leq \frac{1}{2(90-1)} \prod_{i=1}^{88}\frac{1}{(|90-i|-1)^2}\\
<&1.7\times 10^{-271}.
\end{align*}
This shows that all the eigenvalues of $A_2$ that are smaller than $10$ can be hardly perturbed by the subdiagonal $b_{n-k}$ (the same argument shows that more than 80 eigenvalues of $A_2$ can be regarded as converged to within accuracy $10^{-16}$).
\appendix
\section{Multiple eigenvalues}
In the text we assumed that all the eigenvalues of $A+tE$ are simple for all $0\leq t\leq 1$. Here we treat the case where multiple eigenvalues exist, and show that the results we proved hold exactly the same. 

We note that in \cite{barlowdemmel90,barlow00} it is claimed that $A+tE$ can only have multiple eigenvalues on a set of $t$ of measure zero, and hence \eqref{pert} can be integrated on $t$ such that $A+tE$ has only simple eigenvalues. However we cannot use this argument, which can be seen by a simple counterexample $A=E=I$, for which $A+tE$ has a multiple eigenvalue for all $0\leq t\leq 1$. Hence we need a different approach. 
\subsection{Muliple eigenvalue first-order perturbation expansion}
First we review a known result on multiple eigenvalue first-order perturbation expansion \cite{sun92,moro97,holder09}. 
Suppose that a Hermitian matrix $A$ has a multiple eigenvalue $\lambda_0$  of multiplicity $r$, such that there exists a unitary matrix $Q=(Q_1,Q_2)$  such that  
\begin{equation}  \label{xy}
Q^HAQ=\begin{bmatrix} \lambda_0 I&0\\0&\Lambda \end{bmatrix}, 
\end{equation}
where $\Lambda$ is a diagonal matrix that contains eigenvalues that are not $\lambda_0$. 
Then, the matrix  $A+\epsilon E$ has eigenvalues $\widehat\lambda_1,\widehat\lambda_2, \ldots,\widehat\lambda_r$ admitting the first-order expansion
\begin{equation}  \label{expand}
\widehat \lambda_i= \lambda_0 +\mu_i(Q_1^HEQ_1)\epsilon+o(\epsilon), \quad i=1,2,\ldots, r,
\end{equation}
where $\mu_i(Y_1^HEX_1)$ denotes the $i$th eigenvalue of the  $r$-by-$r$ matrix $Q_1^HEQ_1$. 

Using \eqref{expand}, we obtain the partial derivative corresponding to \eqref{pert} when  $A+tE$ has a multiple eigenvalue $\lambda_i(t)=\lambda_{i+1}(t)=\cdots =\lambda_{i+r-1}(t)$  of multiplicity $r$:
\begin{equation}\label{pert2}
\frac{\partial \lambda_{i+j-1}(t)}{\partial t}= \mu_j(Q_1(t)^HE Q_1(t)),\quad 1\leq j\leq r.\end{equation}
Now, let $Q_1(t)^HE Q_1(t)=U^H\Lambda U$ be the eigendecomposition where the diagonals of\\
 $\Lambda=\mbox{diag}(\mu_j(Q_1(t)^HE Q_1(t)))$ are  in descending order. 
Then $\Lambda=UQ_1(t)EQ_1(t)U^H=\tilde Q_1(t)E\tilde Q_1(t)$, where $\tilde Q_1(t)=Q_1(t)U^H$, so  $\mu_j(Q_1(t)^HE Q_1(t))=q_j(t)^HEq_j(t)$, where $q_j(t)$ denotes the $j$th column of $\tilde Q_1(t)$. Now, since any vector of the form $Q_1(t)v$ is an eigenvector corresponding to the eigenvalue $\lambda_i(t)$, so is $q_j(t)$. Hence we can still write the first-order perturbation expansion of each eigenvalue as in \eqref{pert}, and since Lemma \ref{lem} holds regardless of whether $\lambda_i$ is a multiple eigenvalue or not, we conclude that all the results in the text hold exactly the same when multiple eigenvalues exist. 

\subsection{Note on the trailing term}
Here we claim that the expansion can be made sharper in that the trailing term can be $O(\epsilon^2)$ instead of $o(\epsilon)$ as in the known result \eqref{expand}. 
To see this, we write $E=\begin{bmatrix}E_{11}& E_{21}\\E_{21}^H& E_{22} \end{bmatrix}$, and see in \eqref{xy} that 
\[Q^H(A+E)Q=
\begin{bmatrix} \lambda_0 I+Q_1^HE_{11}Q_1&Q_1^HE_{21}Q_2\\Q_2^HE_{21}^HQ_1&\Lambda +Q_2^HE_{22}Q_2 \end{bmatrix}. 
\]
For sufficiently small $E$ there is a positive $gap$ in the spectrums of the matrices $\lambda_0 I+Q_1^HE_{11}Q_1$ and $\Lambda +Q_2^HE_{22}Q_2$. Hence, using the quadratic eigenvalue perturbation  bounds in \cite{rcli05} we see that the $i$th eigenvalue of $Q^H(A+E)Q$ and those of
 $\begin{bmatrix} \lambda_0 I+Q_1^HE_{11}Q_1&0\\0&\Lambda +Q_2^HE_{22}Q_2 \end{bmatrix}$ differ at most by $\displaystyle\frac{2\|E\|_2^2}{gap+\sqrt{gap^2+4\|E\|_2^2}}$. This is of size $O(\epsilon^2)$ because $gap>0$. Therefore we conclude 
\eqref{expand} can be replaced by 
\begin{equation}  \label{expand2}
\widehat \lambda_i= \lambda_0 +\mu_i(Q_1^HEQ_1)\epsilon+O(\epsilon^2), \quad i=1,2,\ldots, r,
\end{equation}

\end{document}